\documentclass[multi]{cambridge7A}
\usepackage[UKenglish]{babel}
\usepackage[longnamesfirst,sectionbib]{natbib}
\usepackage{chapterbib}
\usepackage{amsmath,amssymb,amsthm,amsfonts}
\usepackage{enumerate,url}
\usepackage{calc, xspace}

\setcitestyle{authoryear,round,semicolon}
\allowdisplaybreaks[1]

\theoremstyle{plain}
 \newtheorem{theorem}{Theorem}[section]
\theoremstyle{definition}
 \newtheorem{definition}[theorem]{Definition}

 \renewcommand{\d}{\,\mathrm{d}}

\providecommand*\Index[1]{#1\index{#1}}
\providecommand*\undex[1]{} 

\begin{document}
\alphafootnotes
\author[W.~S.~Kendall]{Wilfrid S.~Kendall\footnotemark }
\chapter[Distribution asymptotics for L\'evy stochastic area couplings]{Coupling time distribution asymptotics for some couplings of the L\'evy stochastic area}
\footnotetext[1]{Department of Statistics, University of Warwick,
  Coventry CV4 7AL; w.s.kendall@warwick.ac.uk,
  \url{http://www.warwick.ac.uk/go/wsk}}
\arabicfootnotes
\contributor{Wilfrid S. Kendall \affiliation{University of Warwick}}
\renewcommand\thesection{\arabic{section}}
\numberwithin{equation}{section}
\renewcommand\theequation{\thesection.\arabic{equation}}

\begin{abstract}
We exhibit some explicit co-adapted couplings for \(n\)-dimensional Brown\-ian motion and all its L\'evy stochastic areas. In the two-dimensional case we show how to derive exact asymptotics for the coupling time under various mixed coupling strategies, using Dufresne's formula for the distribution of exponential functionals of Brownian motion. This yields quantitative asymptotics for the distributions of random times required for
certain simultaneous couplings of stochastic area and Brownian motion. The approach also applies to higher dimensions, but will then lead to upper and lower bounds rather than exact asymptotics.
\end{abstract}

\subparagraph{Keywords}Brownian motion, co-adapted coupling, coupling time distribution, Dufresne formula, exponential functional of Brownian motion, Kolmogorov diffusion, L\'evy stochastic area, maximal coupling, Morse--Thue sequence, non-co-adapted coupling, reflection coupling, rotation coupling, stochastic
   differential, synchronous coupling

\subparagraph{AMS subject classification (MSC2010)}60J65, 60H10

 \section*{Introduction}\label{sec:introduction}\index{Kingman, J. F. C.!influence}
It is a pleasure to present this paper as a homage to my DPhil supervisor John Kingman, in grateful acknowledgement of the formative period which I spent as his research student at Oxford, which launched me into
a deeply satisfying exploration of the world of mathematical research. It seems fitting in this paper to present an overview of 
a particular aspect of probabilistic
\index{coupling|(}coupling theory which has fascinated me for a considerable time; given that one can couple two copies of a random process, when can one in addition couple other associated functionals of the processes? How far can one go?

Motivations for this question include: the sheer intellectual curiosity of discovering just how far one can push the notion of probabilistic coupling; the consideration that coupling has established itself as an extremely powerful tool in probability theory and therefore that any increase in its scope is of potential significance; and the thought that the challenge of coupling in extreme circumstances may produce new paradigms in coupling to complement that of the classic reflection coupling.

It has been known since the 1970s that in principle one can couple two random
processes
\index{coupling!maximal coupling|(}\emph{maximally}; at first encounter this fact continues
to delight and surprise researchers. I summarize this point in
Section \ref{sec:couplings} and also describe the important class of
\index{coupling!co-adapted coupling|(}\emph{co-adapted} couplings. These satisfy more
restrictive requirements than maximal couplings, are typically less efficient,
but are also typically much easier to construct. Since
\index{Lindvall, T.}\citet{Lindvall-1982a}'s seminal preprint we have known
how to couple Euclidean
\index{Brown, R.!Brownian motion|(}Brownian motion using a simple reflection argument in
a way which (most unusually) is both maximal \emph{and} co-adapted, and this
has led to many significant developments and generalizations, some of which
are briefly sketched in Section \ref{sec:reflection}. This leads to
Section \ref{sec:coupling-feature}, which develops the main content of the
paper; what can we now say about the question, how to couple not just Brownian
motion, say, but also associated
\index{path!integral}path integrals? Of course we then need to
vary our strategy, using not just
\index{coupling!reflection coupling}reflection coupling but also so-called
\index{coupling!synchronous coupling}synchronous coupling (in which the two
processes move in parallel), and even
\index{coupling!rotation coupling}rotation coupling, which correlates
different coordinates of the two processes. In \cite{Kendall-2007} I showed
how to couple (co-adaptively) \(n\)-dimensional Brownian motion and all its
\index{Levy, P.@L\'evy, P.!L\'evy stochastic area}stochastic areas; this work is reviewed in Section \ref{subsec:explicit-strategies} using a rather more explicit coupling strategy and then new computations are introduced (in Section \ref{sec:coupling-time}) which establish explicit
asymptotics for the
\index{coupling!time}coupling time for suitable coupling strategies in the two-dimensional case, and which can be used to derive na\"{\i}ve bounds in higher dimensions.
Section \ref{sec:conclusion} concludes the paper with some indications of future research directions.

\section{Different kinds of couplings}\label{sec:couplings}
Probabilistic couplings are used in many different ways: 
couplings (realizations of two random processes on the same probability space) can be constructed so as to arrange
any of the following:
\begin{itemize}
\item  for the two processes to agree after some random time (which is to say,
for the coupling to be \emph{successful}). This follows the pioneering work of
\index{Doeblin, W.}\citet{Doeblin-1938}, which uses this idea to provide a
coupling proof of convergence to equilibrium for finite-state ergodic
\index{Markov, A. A.!Markov chain}Markov chains;
\item  for the two processes to be interrelated by some \emph{monotonicity}
property\break---a common use of coupling in the study of
\index{interacting particle system}interacting particle systems
\index{Liggett, T. M.}\citep{Liggett-2005};
\item  for one process to be linked to the other so as to provide some
informative and fruitful \emph{representation}, as in the case of the
\index{Kingman, J. F. C.!Kingman coalescent}coalescent \citep{Kingman-1982};
\item  for one of the processes to be an illuminating \emph{approximation} to
the other; this appears in an unexpected way in
\index{Barbour, A. D.}\index{Holst, L.}\index{Janson, S.}\citet{BarbourHolstJanson-1992}'s approach to
\index{Stein, C. M.!Stein--Chen approximation}Stein--Chen approximation.
\end{itemize}
These considerations often overlap. Aiming for successful coupling has historical precedence and is in some sense {thematic} for coupling theory, and we will focus on this task here.

\subsection{Maximal couplings}
The first natural question is, how fast can coupling occur? There is a remarkable and satisfying answer, namely that one can in principle construct a coupling which is \emph{maximal} in the sense that it maximizes the probability of coupling before 
time \(t\) for all possible \(t\): see
\index{Griffeath,
D.}\citeauthor{Griffeath-1975} \citetext{\citeyear{Griffeath-1975}, \citeyear{Griffeath-1978}},
\index{Pitman, J. [Pitman, J. W.]}\citet{Pitman-1976},
\index{Goldstein, S.}\citet{Goldstein-1978}. Briefly, maximal couplings
convert the famous
\index{Aldous, D. J.}Aldous inequality (the probability of coupling is bounded above by a simple multiple of the total variation between distributions) into an equality.
Constructions of maximal couplings are typically rather involved, and in
general may be expected to involve demanding potential-theoretic questions
quite as challenging as any problem which the coupling might be supposed to
solve. Pitman's approach is perhaps the most direct, involving explicit
construction of suitable
\index{time reversal}time-reversed Markov chains.

\subsection{Co-adapted coupling}
Maximal couplings being generally hard to construct, it is useful to consider
couplings which are stricter in the sense of requiring the coupled processes
both to be adapted to the same
\index{filtration}filtration. 
Terminology in the literature varies: \emph{Markovian}, when the coupled
processes are jointly Markov, with prescribed marginal kernels
\index{Burdzy, K.}\citep{BurdzyKendall-2000}; \emph{co-immersed}
\index{Emery, M.@\'Emery, M.}\citep{Emery-2005} or \emph{co-adapted} to emphasize the r\^ole of the filtration.
The idea of a co-adapted coupling is simple enough, though its exact
mathematical description is somewhat tedious: here we indicate the definition
for
\index{Markov, A. A.!Markov process}Markov processes.
\begin{definition}\label{def:co-adapted-coupling}
Suppose one is given two continuous-time Markov processes \(X^{(1)}\) and \(X^{(2)}\), with corresponding semigroup kernels defined for bounded measurable functions \(f\) by 
\[
P_t^{(i)}f(z)\quad=\quad\mathbb E\bigl[f(X^{(i)}_{s+t})\bigm|X^{(i)}_s=z,\;X^{(i)}_u \text{ for }u<s\bigr].
\]
A \emph{co-adapted coupling} of \(X^{(1)}\) and \(X^{(2)}\) is a pair of random processes \(\widetilde{X}^{(1)}\) and \(\widetilde{X}^{(2)}\)
defined on the same filtered probability space \((\Omega, \mathfrak{F}, \{\mathfrak{F}_t:t\geq0\}, \mathbb{P})\), 
both adapted to the common filtration \(\{\mathfrak{F}_t:t\geq0\}\) (hence `co-adapted') and satisfying
\[
P_t^{(i)}f(z)\quad=\quad\mathbb E\bigl[f(\widetilde{X}^{(i)}_{s+t})\bigm|\mathfrak{F}_s,\;\widetilde{X}^{(i)}_s=z\bigr]
\]
for \(i=1\), 2, for each bounded measurable function \(f\), each \(z\),
each \(s\), \(t>0\).
\end{definition}
Thus the individual stochastic dynamics of each \(\widetilde{X}^{(i)}\) agree with 
those of the corresponding \({X}^{(i)}\) \emph{even when the past behaviour of the opposite process is also taken into account}. (This is typically not the case for maximal couplings.) In particular, if the \(X^{(i)}\) are Brownian motions then their forward increments are independent of the past given by the filtration.
Moreover if the processes are specified using
\index{stochastic differential equation}stochastic differential equations driven by Brownian motion then general co-adapted couplings can be represented 
using stochastic calculus
(possibly at the price of enriching the filtration), as observed in passing by
\index{Emery, M.@\'Emery, M.}\citet{Emery-2005}, and as described more
formally in \citet[Lemma 6]{Kendall-2009a}. Briefly, any co-adapted coupling
of vector-valued Brownian motions \(\underline{A}\) and \(\underline{B}\) can
be represented by expressing \(\underline{A}\) as a
\Index{stochastic integral} with respect to \(\underline{B}\) and perhaps another independent Brownian motion  \(\underline{C}\): we use this later at Equation \eqref{eq:coupling-equation}.

\subsection{Coupling at different times}
There are many other useful couplings falling outside this framework: for
example
\index{Thorisson, H.}\citet{Thorisson-1994} discusses the idea of a
\index{coupling!shift-coupling}\emph{shift-coupling}, which weakens the coupling requirement by permitting processes to couple at \emph{different} times;
\citet{Kendall-1994} uses co-adapted coupling of \emph{time-changed} processes
as part of an exploration of regularity for
\index{harmonic map}harmonic maps. However in this paper we will focus on co-adapted couplings.

\section{Reflection coupling}\label{sec:reflection}\index{coupling!reflection coupling|(}
The dominant example of coupling is \emph{reflection coupling} for Euclidean Brownian motions \(\underline{A}\) and \(\underline{B}\),
dating back to
\index{Lindvall, T.|(}\citet{Lindvall-1982a}'s preprint:
construct \(\underline{A}\) from \(\underline{B}\) by
reflecting \(\underline{B}\) in the line segment running
from \(\underline{B}\) to \(\underline{A}\). That this is a maximal coupling
follows from an easy computation involving the
\Index{reflection principle}.
It can be expressed as a co-adapted coupling; 
 the Brownian increment for \(\underline{A}\) is derived from that of \(\underline{B}\) by a reflection in the line segment running from \(\underline{B}\) to \(\underline{A}\). 
Many modifications of the reflection coupling have been derived to cover various situations; we provide a quick survey in the remainder of this section.

\subsection{Maximality and non-maximality}
The reflection coupling is unusual in being both co-adapted and maximal.
\index{Hsu, E. P.}\index{Sturm, K.-T.}\citet{HsuSturm-2003} point out that
reflection coupling fails to be maximal even for  Euclidean Brownian motion if
the Brownian motion is stopped on exit from a prescribed domain.
\index{Kuwada, K.}(\citealt{KuwadaSturm-2007} discuss the
\index{manifold|(}manifold case; see also \citealt{Kuwada-2007}, \citeyear{Kuwada-2009}.)
Perhaps the simplest example of an instance where no co-adapted coupling can
be maximal arises in the case of the
\index{Ornstein, L. S.!Ornstein--Uhlenbeck process}Ornstein--Uhlenbeck process
\index{Connor, S. B.}\citep[PhD thesis, Theorem 3.15]{Connor-2007a}. Consider the problem of constructing successful co-adapted couplings between (i) an Ornstein--Uhlenbeck process begun at \(0\), and (ii) an Ornstein--Uhlenbeck process run in statistical equilibrium. A direct argument shows that no such co-adapted coupling can be maximal; however in this case reflection coupling is maximal amongst all \emph{co-adapted} couplings.
(The study of couplings which are
\index{coupling!maximal coupling|)}maximal in the class of co-adapted couplings
promises to be an interesting field: the case of
\index{random walk (RW)}random walk on the \Index{hypercube} is treated by
\index{Jacka, S. D.}\citealt{ConnorJacka-2008}.)

\subsection{Coupling for Diffusions}
A variant of reflection coupling for
\index{diffusion!elliptic diffusion}elliptic diffusions with smooth coefficients is
discussed in
\index{Lindvall, T.|)}\index{Rogers, L. C. G.}\citet{LindvallRogers-1986}
and further in
\index{Chen, M. F.}\index{LiSF@Li, S. F.}\citet{ChenLi-1989}. `Reflection' here depends on interaction between the two diffusion matrices, and in general the two coupled diffusions do not play symmetrical r\^oles.
In the case of Brownian motion on a
\index{manifold|)}manifold one can use the mechanisms of stochastic
development and stochastic
\Index{parallel transport} to define co-adapted couplings in a more symmetrical manner.
The behaviour of general co-adapted Brownian couplings on Riemannian manifolds
is related to curvature. \citet{Kendall-1986a} shows that successful
co-adapted coupling can never be almost-surely successful in the case of a
simply-connected manifold with negative
\index{curvature}curvatures bounded away from zero. On the other hand a
geometric variant of reflection coupling known as
\index{coupling!mirror coupling}\emph{mirror coupling} will always be almost-surely successful if the manifold has non-negative Ricci curvatures
(\citealt{Kendall-1986b},
\index{Cranston, M.}\citealt{Cranston-1992}, \citealt{Kendall-1988d}).
Indeed
\index{Renesse, M. von}\Citet{Renesse-2004} shows how to generalize mirror coupling even to non-manifold contexts.

\section{Coupling more than one feature of the process}\label{sec:coupling-feature}
The particular focus of this paper concerns ongoing work on the following question: is it possible co-adaptively to couple more than one feature of a random process at once? To be explicit, is it possible to couple not just the location but also some functional of the path?

On the face of it, this presents an intimidating challenge: control of
difference of
\index{path!functional}path functionals 
by coupling is necessarily indirect and it is possible that all attempts to control the discrepancy between path functionals will inevitably jeopardize coupling of the process itself.

Further thought shows that it is sensible to confine attention to cases where
the process together with its functional form a
\index{diffusion!hypoelliptic diffusion}\emph{hypoelliptic diffusion}, since
in such cases the
\index{Hormander, L.@H\"ormander, L.!H\"ormander regularity theorem}H\"ormander regularity theorem guarantees existence of a
density, and this in turn shows that general coupling (not necessarily
non-co-adapted) is possible in principle.
\index{Hairer, M.}(\citealt{Hairer-2002} uses this approach to produce
non-co-adapted couplings, using careful analytic estimates and a regeneration
argument which corresponds to \citealt{Lindvall-2002}'s
\index{coupling!gamma coupling@$\gamma$-coupling}`\(\gamma\) coupling'.)

Furthermore it is then natural to restrict attention to diffusions with nilpotent group symmetries, where one may hope most easily to discover co-adapted couplings which will be susceptible to extensive generalization, parallelling the generalizations of the Euclidean reflection coupling which have been described briefly in Section \ref{sec:reflection}.

\subsection{Kolmogorov diffusion}\index{Kolmogorov, A. N.!Kolmogorov diffusion}

Consider the so-called Kolmogorov diffusion: scalar Brownian motion \(B\) plus
its time integral \(\int B\d t\). This determines a simple
\index{diffusion!nilpotent diffusion|(}nilpotent diffusion, and in fact it can be coupled co-adaptively by varying sequentially
between reflection coupling and
\index{coupling!synchronous coupling|(}\emph{synchronous} coupling (allowing
the two Brownian motions to move in parallel) as shown
in
\index{Ben Arous, G.|(}\index{Cranston, M.}\citet{BenArousCranstonKendall-1995}.
\index{Jansons, K. M.}\index{Metcalfe, P. D.}\citet{JansonsMetcalfe-2005b}
describe some numerical investigations concerned with optimizing an
exponential moment of the
\index{coupling!time}coupling time.

The idea underlying this coupling is rather simple. Suppose that we wish to
couple \((B, \int B\d t)\) with \((\widetilde{B},\int\widetilde{B}\d
t)\). Set \(U=\widetilde{B}-B\) and \(V=\int\widetilde{B}\d t-\int B\d
t\). Co-adapted couplings include
\Index{stochastic integral} representations such as \(\d\widetilde{B}=J\d B\),
for co-adapted \(J\in\{-1,1\}\); \(J=1\) yields synchronous coupling
and \(J=-1\) yields reflected coupling. Suppose \(U_0\neq0\) and \(V=0\)
(which can always be achieved by starting with a little reflected or
synchronous coupling unless \(U=V=0\) from the start, in which case nothing
needs to be done). We can cause \((U,V)\) to wind repeatedly around \((0,0)\)
in ever smaller loops as follows: first use reflection coupling till \(U\)
hits \(-U_0/2\), then synchronous coupling till \(V\) hits \(0\), then repeat
the procedure but substituting \(-U_0/2\) for \(U_0\). A
\undex{Borel, F. E. J. E.!Borel--Cantelli lemma}Borel--Cantelli argument combined with Brownian scaling shows that \((U,V)\) then hits \((0,0)\) in finite time.

\index{Price, C. J.}\citet{KendallPrice-2004} present a cleaned-up version of this argument (together with an extension to deal in addition with \(\int\!\!\int B\d s\,\d t\)). 

Curiously, this apparently artificial example can actually be applied to the
study of the
\index{tail sigma algebra@tail $\sigma$-algebra}tail \(\sigma\)-algebra of a
certain
\index{diffusion!relativistic diffusion}relativistic diffusion discussed by
\index{Bailleul, I.}\citet{Bailleul-2008}.

Remarkably it is possible to do very much better by using a completely
different and implicit method: one can couple not just the time integral, but
also any finite number of additional iterated time
integrals \citep{KendallPrice-2004}, by concatenating reflection and
synchronous couplings using the celebrated
\index{Morse, H. C.!Morse--Thue binary sequence}\emph{Morse--Thue} binary
sequence \(0110100110010110\ldots\). Scaled iterations of state-dependent
perturbations of the resulting concatenation of couplings can be used to
deliver coupling in a finite time; the perturbed Morse--Thue sequences encode indirect controls of higher-order iterated integrals. 

\subsection{L\'evy stochastic areas}\index{Levy, P.@L\'evy, P.!L\'evy stochastic area|(}
Moving from scalar to planar Brownian motion, the natural question is now
whether one can co-adaptively couple the
\index{diffusion!nilpotent diffusion|)}nilpotent diffusion formed by Brownian
motion \((B_1, B_2)\) and the L\'evy stochastic area \(\int(B_1\d B_2-B_2\d
B_1)\). This corresponds to coupling a
\index{diffusion!hypoelliptic diffusion}hypoelliptic Brownian motion on the
\index{Heisenberg, W. K.!Heisenberg group}Heisenberg group, and 
\index{Ben Arous, G.|)}\index{Cranston, M.}\citeauthor{BenArousCranstonKendall-1995} determine an explicit successful
coupling based on extensive explorations using
\Index{computer algebra}.

Again one can do better \citep{Kendall-2007}. Not only can one construct a
simplified coupling for the \(2\)-dimensional case based only on reflection
and synchronous couplings (switching from reflection to synchronous coupling
whenever a geometric difference of the stochastic areas exceeds a fixed
multiple of the squared distance between the two coupled Brownian motions),
but also one can successfully couple \(n\)-dimensional Brownian motion plus
a \(\tbinom{n}{2}\)  basis of the various stochastic areas. In the remainder
of this paper we will indicate the method used, which moves beyond the use of
\index{coupling!reflection coupling|)}reflection and
\index{coupling!synchronous coupling|)}synchronous couplings to involve
\index{coupling!rotation coupling}\emph{rotation couplings} as well (in which coordinates of one of the Brownian motions can be correlated to quite different coordinates of the other).

\subsection{Explicit strategies for coupling L\'evy stochastic area}\label{subsec:explicit-strategies}
Here we describe a variant coupling strategy for the \(n\)-dimensional case  which is more explicit than the strategy proposed in \citet{Kendall-2007}.
As described in \citet[Lemma 6]{Kendall-2007}, suppose that \(\underline{A}\)
and \(\underline{B}\) are co-adaptively coupled \(n\)-dimensional Brownian
motions. Arguing as in \citet{Kendall-2009a}, and enriching the
\index{filtration}filtration if necessary, we may represent any such coupling in terms of a further \(n\)-dimensional Brownian motion \(\underline{C}\), independent of \(\underline{B}\);
\begin{equation}\label{eq:coupling-equation}
 \d\underline{A} \quad=\quad 
 \underline{\underline{J}}^\top\d\underline{B} + \underline{\underline{\widetilde{J}}}^\top\d\underline{C}\,,
\end{equation} 
where \(\underline{\underline{J}}\), \(\underline{\underline{\widetilde{J}}}\) are predictable \((n\times n)\)-matrix-valued processes satisfying the constraint
\begin{equation}
\underline{\underline{J}}^\top\,\underline{\underline{J}}\;+\;\underline{\underline{\widetilde{J}}}^\top\,\underline{\underline{\widetilde{J}}}
\quad=\quad
\underline{\underline{\mathbb{I}}}\,
\label{eq:brownian-relation}
\end{equation} 
and where \(\underline{\underline{\mathbb{I}}}\) represents the \((n\times n)\) identity matrix.

Note that the condition \eqref{eq:brownian-relation} is equivalent to the following set of symmetric matrix inequalities
for the
\index{coupling!co-adapted coupling|)}co-adapted process \(\underline{\underline{J}}\) (interpreted in a spectral sense):
\begin{equation}
\underline{\underline{J}}^\top\,\underline{\underline{J}}
\quad\leq\quad
\underline{\underline{\mathbb{I}}}\,.
\label{eq:coupling-condition}
\end{equation}

Thus a legitimate
\index{coupling!control|(}\emph{coupling control} \(\underline{\underline{J}}\) must take values in a compact convex set of \(n\times n\) matrices defined by \eqref{eq:coupling-condition}.

The matrix process \(\underline{\underline{J}}\) measures the
correlation \((\d\underline{B}\,\d\underline{A}^\top)/\d t\)\break bet\-ween
the
\index{Brown, R.!Brownian differential}Brownian differentials \(\d\underline{A}\) and \(\d\underline{B}\); for
convenience
let \(\underline{\underline{S}}=\tfrac{1}{2}(\underline{\underline{J}}+\underline{\underline{J}}^\top)\)
and \(\underline{\underline{A}}=\tfrac{1}{2}(\underline{\underline{J}}-\underline{\underline{J}}^\top)\)
refer to the symmetric and skew-symmetric parts
of \(\underline{\underline{J}}\). The coupling problem solved
in \citet{Kendall-2007} is to choose an
adapted \(\underline{\underline{J}}=\underline{\underline{S}}+\underline{\underline{A}}\)
which brings \(\underline{A}\) and \(\underline{B}\) into agreement at a
\index{coupling!time}coupling time \(T_\text{coupling}\) which is simultaneously a coupling time for all the \(\tbinom{n}{2}\) corresponding pairs of stochastic area integrals \(\int\left(A_i\d A_j - A_j\d A_i\right)\) and \(\int\left(B_i\d B_j - B_j\d B_i\right)\).

To measure progress towards this simultaneous coupling, set \(\underline{X}=\underline{A}-\underline{B}\) 
and define \(\underline{\underline{\mathfrak{A}}}\) to be a skew-symmetric matrix of geometric differences between stochastic areas with \((i,j)^\text{th}\) entry
\begin{equation}\label{eq:areal-difference}
\mathfrak{A}_{ij} =
\int\left(A_i\d A_j - A_j\d A_i\right)
-
\int\left(B_i\d B_j - B_j\d B_i\right)
+A_i B_j - A_j B_i.
\end{equation} 

The nonlinear correction term \(A_i B_j - A_j B_i\) is important because it converts \(\mathfrak{A}_{ij}\) into a geometrically natural quantity, invariant under shifts of coordinate system, and also because it supplies a very useful contribution to the drift in the subsequent It\^o analysis.
Of course \(\underline{\underline{\mathfrak{A}}}\) and \(\underline{X}\) both vanish at a given time \(t\) if and only if at that time both \(\underline{A}=\underline{B}\) (so in particular the correction term vanishes) and also all the corresponding stochastic areas agree.

Some detailed
\index{Ito, K.@It\^o, K.!It\^o calculus}It\^o calculus (originally carried out
in an implementation of the It\^o calculus
\index{computer algebra!Axiom@{\it Axiom}} procedures \emph{Itovsn3} in
\emph{Axiom}, \citealp{Kendall-2001b}, but now checked comprehensively  by
hand) can now be used to derive the following system of
\index{stochastic differential equation}stochastic differential equations for the squared distance \(V^2=\|\underline{X}\|^2\) and the `squared areal difference' \(U^2=\operatorname{trace}(\underline{\underline{\mathfrak{A}}}^\top\underline{\underline{\mathfrak{A}}})={\sum\sum}_{ij} \mathfrak{A}_{ij}^2\):
\begin{align}
(\d (V^2))^2 \quad&=\quad 8 \,\underline{\nu}^\top\left(\underline{\underline{\mathbb{I}}}-\underline{\underline{S}}\right)\underline{\nu} \; V^2\d t\,,
\nonumber\\
\operatorname{Drift}\d (V^2) \quad&=\quad 2 \operatorname{trace}\left(\underline{\underline{\mathbb{I}}}-\underline{\underline{S}}\right) \; \d t\,,
\nonumber\\
\d (V^2)\times\d(U^2) \quad&=\quad -16\, \underline{\nu}^\top\underline{\underline{Z}}^\top\underline{\underline{A}}\,\underline{\nu} \; U V^2\d t\,,
\nonumber\\
(\d (U^2))^2 \quad&=\quad 32\, \underline{\nu}^\top\underline{\underline{Z}}^\top\left(\underline{\underline{\mathbb{I}}}+\underline{\underline{S}}\right)\underline{\underline{Z}}\,\underline{\nu} \; U^2 V^2\d t\,,
\nonumber\\
\operatorname{Drift}\d (U^2) \quad&=\quad 4 \operatorname{trace}\left(\underline{\underline{Z}}^\top\underline{\underline{A}}\right) \; U\d t
\;+\;
\nonumber\\
& \quad\qquad\quad
\;+\;
4\left(
\operatorname{trace}\left(\underline{\underline{\mathbb{I}}}+\underline{\underline{S}}\right)
-
\underline{\nu}^\top\left(\underline{\underline{\mathbb{I}}}+\underline{\underline{S}}\right)\underline{\nu}
\right)\; V^2\d t\,.
\label{eq:sde-V2-U2}
\end{align}
Here the vector \(\underline{\nu}\) and the
matrix \(\underline{\underline{Z}}\) encode relevant underlying geometry:
respectively \(V \underline{\nu}=\underline{X}\)
and \(U\underline{\underline{Z}}=\underline{\underline{\mathfrak{A}}}\). Note
that \(\underline{\nu}\) is a unit vector and \(\underline{\underline{Z}}\)
has unit
\index{Hilbert, D.!Hilbert--Schmidt norm}Hilbert--Schmidt norm: \(\operatorname{trace}{\underline{\underline{Z}}^\top\underline{\underline{Z}}=1}\).

The strategy is to consider \(V^2\) and \(U^2\) on a log-scale: further
stochastic calculus together with suitable choice of coupling
control \(\underline{\underline{J}}\) then permits comparison to two Brownian
motions with constant negative drift in a new time-scale, and a
stochastic-calculus argument shows that \(K=\tfrac{1}{2}\log(V^2)\)
and \(H=\tfrac{1}{2}\log(U^2)\) reach \(-\infty\) at finite time in the
original time-scale measured by \(t\). In fact further stochastic calculus,
based on the
\index{martingale!martingale differential identity}martingale differential identity
\begin{equation*}
 \d \log{Z} \quad=\quad \frac{\d Z}{Z} - \tfrac{1}{2} \left(\frac{\d Z}{Z}\right)^2\,,
\end{equation*}
shows that
\begin{align}
 (\d K)^2 \quad&=\quad 
\tfrac{1}{2}\;\underline{\nu}^\top\left(\underline{\underline{\mathbb{I}}}-\underline{\underline{S}}\right)\underline{\nu} \;\d\widetilde{\tau}\,,
\nonumber\\
 \operatorname{Drift}(\d K) \quad&=\quad 
\tfrac{1}{4}\left(
\operatorname{trace}\left(\underline{\underline{\mathbb{I}}}-\underline{\underline{S}}\right)
-2\;\underline{\nu}^\top\left(\underline{\underline{\mathbb{I}}}-\underline{\underline{S}}\right)\underline{\nu}\right) \;\d\widetilde{\tau}\,,
\nonumber\\
 \d K\times\d H \quad&=\quad 
-\;\underline{\nu}^\top\underline{\underline{Z}}^\top\underline{\underline{A}}\,\underline{\nu} \;\tfrac{\d\widetilde{\tau}}{W}\,,
\nonumber\\
 (\d H)^2 \quad&=\quad 
2\;\underline{\nu}^\top\underline{\underline{Z}}^\top\left(\underline{\underline{\mathbb{I}}}+\underline{\underline{S}}\right)\underline{\underline{Z}}\,\underline{\nu} \; \tfrac{\d\widetilde{\tau}}{W^2}\,,
\nonumber\\
 \operatorname{Drift}(\d H) \quad&=\quad 
\tfrac{1}{2}\;\operatorname{trace}\left(\underline{\underline{Z}}^\top\underline{\underline{A}}\right) \; \tfrac{\d\widetilde{\tau}}{W}
\;+\;
\nonumber\\
\;+\;
\tfrac{1}{2}&
\left(
\operatorname{trace}\left(\underline{\underline{\mathbb{I}}}+\underline{\underline{S}}\right)
-
\underline{\nu}^\top\left(\underline{\underline{\mathbb{I}}}+\underline{\underline{S}}\right)\underline{\nu}
-4\;\underline{\nu}^\top\underline{\underline{Z}}^\top\left(\underline{\underline{\mathbb{I}}}+\underline{\underline{S}}\right)\underline{\underline{Z}}\,\underline{\nu}
\right)\;
\tfrac{\d\widetilde{\tau}}{W^2}\,.\label{eq:sde-K-H}
\end{align}
Here \(\d\widetilde{\tau}=4\tfrac{\d t}{V^2}\) determines the new time-scale \(\widetilde{\tau}\), and \(W=U/V^2=\exp\left(H-2K\right)\). It is clear from the system \eqref{eq:sde-K-H} that the contribution of \(\tfrac{1}{2}\,\operatorname{trace}\bigl(\underline{\underline{Z}}^\top\underline{\underline{A}}\bigr)\; \tfrac{\d\widetilde{\tau}}{W}\) (deriving ultimately from the areal difference correction term in \eqref{eq:areal-difference}) is potentially a flexible and effective component of the control if \(1/W^2\) is small. On the other hand if this is to be useful then \(\underline{\underline{\mathbb{I}}}-\underline{\underline{S}}\) must be correspondingly reduced so as to ensure that \(H\) and \(K\) are subject to dynamics on comparable time-scales.

\citet{Kendall-2007} considers the effect of a control \(\underline{\underline{J}}\) which is an affine mixture of \emph{reflection} \(\underline{\mathbb{I}}-2\,\underline{\nu}\,\underline{\nu}^\top\) and \emph{rotation} \(\exp\left(-\theta\, \underline{\underline{Z}}\right)\). Second-order Taylor series expansion of the matrix exponential is used to overcome analytical complexities at the price of some mild asymptotic analysis. Here we indicate an alternative route,
replacing the second-order truncated expansion
\(\underline{\underline{J}}^\prime=\underline{\underline{\mathbb{I}}}-\theta\, \underline{\underline{Z}} -\tfrac{1}{2}\theta^2 \,\underline{\underline{Z}}^\top\underline{\underline{Z}}\) (which itself fails to satisfy \eqref{eq:coupling-condition} and is not therefore a valid coupling control) by 
\(\underline{\underline{J}}^{\prime\prime}=\underline{\mathbb{I}}-\theta\, \underline{\underline{Z}} -\theta^2 \,\underline{\underline{Z}}^\top\underline{\underline{Z}}\), which does satisfy \eqref{eq:coupling-condition} when \(\theta^2\leq2\) (using the fact that non-zero eigenvalues of \(\underline{\underline{Z}}^\top\underline{\underline{Z}}\) have multiplicity \(2\), so that \(\operatorname{Trace}(\underline{\underline{Z}}^\top\underline{\underline{Z}})=1\) implies that \(\underline{\underline{0}}\leq\underline{\underline{Z}}^\top\underline{\underline{Z}}\leq\tfrac{1}{2}\,\underline{\underline{\mathbb{I}}}\)). Thus we can consider
\begin{equation}
 \underline{\underline{J}}\quad=\quad
\underline{\underline{\mathbb{I}}} - 2 p\,\underline{\nu}\,\underline{\nu}^\top -(1-p)\theta\,\underline{\underline{Z}}
-(1-p)\theta^2\,\underline{\underline{Z}}^\top\underline{\underline{Z}}\,,
\end{equation}
which is a valid coupling control satisfying \eqref{eq:coupling-condition} when \(0\leq p\leq1\) and \(\theta^2\leq2\).

This leads to the following
\index{stochastic differential system|(}stochastic differential system,
where terms which will eventually be negligible have been separated out:
\begin{align}
 (\d K)^2 \quad&=\quad \left(p+\tfrac{1}{2}\left\|\underline{\underline{Z}}\,\underline{\nu}\right\|^2\theta^2\right) \; \d\widetilde{\tau}
\;-\;
\tfrac{p}{2}\left\|\underline{\underline{Z}}\,\underline{\nu}\right\|^2\theta^2
\;\d\widetilde{\tau}\,,
&\nonumber\\
 \operatorname{Drift}(\d K) \quad&=\quad -\tfrac{1}{2} \left(p-\tfrac{1}{2}\left(1-2\|\underline{\underline{Z}}\,\underline{\nu}\|^2\right)\theta^2\right) \d\widetilde{\tau}
\nonumber \\
& \;-\;
\tfrac{p}{4}\left(1-2\|\underline{\underline{Z}}\,\underline{\nu}\|^2\right)\theta^2
 \;\d\widetilde{\tau}\,,
&\nonumber\\
 \d K\times\d H \quad&=\quad \|\underline{\underline{Z}}\,\underline{\nu}\|^2\theta
\;\tfrac{\d\widetilde{\tau}}{W}
\;-\:
p\|\underline{\underline{Z}}\,\underline{\nu}\|^2\theta
\;\tfrac{\d\widetilde{\tau}}{W}\,,
&\nonumber\\
 (\d H)^2 \quad&=\quad 
4\;\|\underline{\underline{Z}}\,\underline{\nu}\|^2\; \tfrac{\d\widetilde{\tau}}{W^2}
\;-\; 2({1-p})\|\underline{\underline{Z}}^\top\underline{\underline{Z}}\,\underline{\nu}\|^2\theta^2
\; \tfrac{\d\widetilde{\tau}}{W^2}\,,
&\nonumber\\
 \operatorname{Drift}(\d H) \quad&=\quad 
-\tfrac{1}{2}\theta
\; \tfrac{\d\widetilde{\tau}}{W} 
\;+\; \left(n-1-4\|\underline{\underline{Z}}\,\underline{\nu}\|^2\right)
\;\tfrac{\d\widetilde{\tau}}{W^2}
&\nonumber\\
&\;+\;\tfrac{1}{2}p\theta\; \tfrac{\d\widetilde{\tau}}{W}
\;-\;\tfrac{1-p}{2}\left(
1
-\left\|\underline{\underline{Z}}\,\underline{\nu}\right\|^2
-4\;
\|\underline{\underline{Z}}^\top\underline{\underline{Z}}\,\underline{\nu}\|^2
\right)\theta^2\;
\tfrac{\d\widetilde{\tau}}{W^2}\,.
&\label{eq:sde-K-H-realized}
\end{align}
In order to ensure comparable dynamics for \(K\) and \(H\), set \(p=\alpha^2/W^2\) and \(\theta=\beta/W\) (valid when \(W\geq\max\bigl\{\tfrac{1}{\sqrt{2}}\beta,\alpha\bigr\}\)): writing \(\d\tau=\d\widetilde{\tau}/W^2=4(V/U)^2\d t\) leads to
\begin{align}
 (\d K)^2 \quad&=\quad \left(\alpha^2+\tfrac{1}{2}\left\|\underline{\underline{Z}}\,\underline{\nu}\right\|^2\beta^2
\;-\;\tfrac{1}{2}\left\|\underline{\underline{Z}}\,\underline{\nu}\right\|^2\tfrac{\alpha^2\beta^2}{W^2}
\right) \; \d{\tau}\,,
\nonumber\\
 \operatorname{Drift}(\d K) \quad&=\quad
 -\mu_1 \d{\tau}\,,
\nonumber\\
 \d K\times\d H \quad&=\quad \|\underline{\underline{Z}}\,\underline{\nu}\|^2\left(1-\tfrac{\alpha^2}{W^2}\right)\beta
\;\d{\tau}\,,
\nonumber\\
 (\d H)^2 \quad&=\quad
4
\;
\left(\|\underline{\underline{Z}}\,\underline{\nu}\|^2 -
 \tfrac{1}{2}({1-\tfrac{\alpha^2}{W^2}})\|\underline{\underline{Z}}^\top\underline{\underline{Z}}\,\underline{\nu}\|^2\tfrac{\beta^2}{W^2}\right)\; {\d{\tau}}\,,
\nonumber\\
 \operatorname{Drift}(\d H) \quad&=\quad 
-\mu_2
\; {\d{\tau}}\,,
\label{eq:sde-K-H-realized2}
\end{align}
where \(\mu_1\) and \(\mu_2\) are given by
\begin{align}
\mu_1 \quad&=\quad \tfrac{1}{2} \alpha^2-\tfrac{1}{4}\left(1-2\|\underline{\underline{Z}}\,\underline{\nu}\|^2\right)\beta^2 
+
\tfrac{1}{4 W^2}\left(1-2\|\underline{\underline{Z}}\,\underline{\nu}\|^2\right)
\alpha^2\beta^2,
\nonumber\\
\mu_2 \quad&=\quad \tfrac{1}{{2}}\beta-n+1+4\|\underline{\underline{Z}}\,\underline{\nu}\|^2 
\nonumber\\ & \qquad
-
\tfrac{1}{2W^2}\left(
\alpha^2\beta\; 
\;-\;\left(1-\tfrac{\alpha^2}{W^2}\right)
\left(
 1
 -\left\|\underline{\underline{Z}}\,\underline{\nu}\right\|^2
 -4\;
 \|\underline{\underline{Z}}^\top\underline{\underline{Z}}\,\underline{\nu}\|^2
 \right)\beta^2
\right)\,.
\label{eq:sde-K-H-realized2-parameters}
\end{align}

In order to fulfil the underlying strategy, \(\mu_1\) and \(\mu_2\) should be chosen to accomplish the following:
\begin{enumerate}[1.]
 \item \(W=\exp(H-2K)\) must remain large; this follows by a
\index{strong law of large numbers (SLLN)}strong-law-of-large-numbers argument if \(\mu_1\), \(\mu_2\) are chosen so that
\(2\mu_1-\mu_2=\operatorname{Drift}\d(H-2K)/\d\tau\) is positive and bounded away from zero
and yet \((\d H)^2/\d\tau\) and \((\d K)^2/\d\tau\) are bounded.
\item Both of \(-\mu_1=\operatorname{Drift}\d(K)/\d\tau\) and
\(-\mu_2=\operatorname{Drift}\d(H)/\d\tau\) must remain negative and bounded away from zero, so that \(K\) and \(H\) both tend to \(-\infty\) as \(\tau\to\infty\).
\item If coupling is to happen in finite time on the \(t\)-time-scale then
\begin{multline*}
T_\text{coupling}\quad=\quad \frac{1}{4}\int_0^\infty \left(\frac{U}{V}\right)^2\d\tau 
\\
\quad=\quad
\frac{1}{4}\int_0^\infty \exp\left(2H - 2K\right)\d\tau 
\quad<\quad \infty\,.
\end{multline*}
This follows almost surely if \(\mu_1\), \(\mu_2\) are chosen so that \(2\mu_1-2\mu_2=\operatorname{Drift}\d(2H-2K)/\d\tau\) is negative and bounded away from zero.
\end{enumerate}
Now the
\Index{inverse function theorem} can be applied to show that for any constant \(L>0\) there is a constant \(L^\prime>0\) such that if \(W^2\geq L^\prime\) then \eqref{eq:sde-K-H-realized2-parameters}
can be solved for
any prescribed \(0\leq\mu_1\), \(\mu_2\leq L\) using \(\alpha^2\), \(\beta^2\leq L^\prime\) (incidentally thus bounding \((\d H)^2/\d\tau\) and \((\d K)^2/\d\tau\)). Moreover by choosing \(L^\prime\) large enough
it then follows that \(p=\alpha^2/W^2\leq1\) and \(\beta^2/W^2\leq 2\), so that the desired \(\mu_1\), \(\mu_2\) can be attained using a valid coupling control.

A comparison with Brownian motion of constant drift now shows that if initially \(W\geq 2L^\prime\) then there is a positive chance that \(W\geq L^\prime\) for all \(\tau\), and thus that coupling happens at \(\tau\)-time infinity, and actual time \(T_\text{coupling}<\infty\). Should \(W\) drop below \(L^\prime\), then one can switch to the pure reflection control (\(p=1\)) and run this control until \(W=U/V^2\) rises again to level \(2L^\prime\). This is almost sure to happen eventually, since otherwise \(V\to0\) and thus \(U=WV^2\to0\), which can be shown to have probability \(0\) under this control.

\subsection{Estimates for coupling time distribution}\label{sec:coupling-time}\index{coupling!time|(}
In the planar case \(n=2\) the stochastic differential system
(\ref{eq:sde-K-H}) under mixed rotation-reflection controls can be simplified
substantially. In this section we go beyond the work of
\index{Ben Arous, G.}\index{Cranston, M.}\citeauthor{BenArousCranstonKendall-1995}(\citeyear{BenArousCranstonKendall-1995}) and \citet{Kendall-2007} by using this simplification to identify limiting distributions for suitable scalings of the coupling time \(T_\text{coupling}\). The
simplification arises because the non-zero eigenvalues of the skew-symmetric matrix \(\underline{\underline{Z}}\) have even multiplicity; consequently in the two-dimensional case we may deduce from \(\operatorname{trace}(\underline{\underline{Z}}^\top\underline{\underline{Z}})=1\) and \(\|\underline{\nu}\|=1\) that
\(\left\|\underline{\underline{Z}}\,\underline{\nu}\right\|^2=\tfrac{1}{2}\) and
\(\|\underline{\underline{Z}}^\top\underline{\underline{Z}}\,\underline{\nu}\|^2=\tfrac{1}{4}\). 
Moreover an
\index{Euler, L.!Euler formula}Euler formula follows from \(\underline{\underline{Z}}\,\underline{\underline{Z}}=-\tfrac{1}{2}\underline{\underline{\mathbb{I}}}\):
\[
\exp\left(-\sqrt{2}\,\theta\,\underline{\underline{Z}}\right)\quad=\quad
\cos\theta\,\underline{\underline{\mathbb{I}}}
-\sqrt{2}\,\sin\theta\,\underline{\underline{Z}}\,.
\]
Accordingly in the planar case the mixed coupling control
\begin{multline}
  \underline{\underline{J}}\quad=\quad
p\left(\underline{\underline{\mathbb{I}}} - 2 \,\underline{\nu}\,\underline{\nu}^\top\right) +q\exp\left(-\sqrt{2}\,\theta\,\underline{\underline{Z}}\right)
 \\
\quad=\quad
\underline{\underline{\mathbb{I}}}
- 2 p \,\underline{\nu}\,\underline{\nu}^\top
- q\,\left(1-\cos\theta\right)\underline{\underline{\mathbb{I}}}
- \sqrt{2}\,q\,\sin\theta\, \underline{\underline{Z}}
\label{eq:mixed-control}
\end{multline} 
(for \(0\leq p=1-q\leq 1\), unrestricted \(\theta\)) renders the system (\ref{eq:sde-K-H}) as
\begin{align}
 (\d K)^2 \quad&=\quad 
(p+\tfrac{1}{2}q(1-\cos\theta) \;\d\widetilde{\tau}\,,
\nonumber\\
 \operatorname{Drift}(\d K) \quad&=\quad 
-\tfrac{1}{2}p \;\d\widetilde{\tau}\,,
\nonumber\\
 \d K\times\d H \quad&=\quad 
\tfrac{1}{\sqrt{2}} q \sin\theta \;\tfrac{\d\widetilde{\tau}}{W}\,,
\nonumber\\
 (\d H)^2 \quad&=\quad 
(2-q(1-\cos\theta)) \; \tfrac{\d\widetilde{\tau}}{W^2}\,,
\nonumber\\
 \operatorname{Drift}(\d H) \quad&=\quad 
-\tfrac{1}{\sqrt{2}}q\sin\theta\;\tfrac{\d\widetilde{\tau}}{W}
-(1-\tfrac{1}{2}q(1-\cos\theta))
\;\tfrac{\d\widetilde{\tau}}{W^2}\,.\label{eq:sde-K-H-planar-case}
\end{align}
If we set \(\theta=0\) (so that the
\index{coupling!control|)}coupling control is a mixture of
\index{coupling!reflection coupling|(}reflection coupling and
\index{coupling!synchronous coupling|(}\emph{synchronous} couplings) then the result can be made to yield a successful coupling:
choosing \(p=\min\{1, \alpha^2/W^2\}\) the
\index{stochastic differential system|)}stochastic differential system becomes
\begin{align}
 (\d K)^2 \quad&=\quad 
\min\{W^2, \alpha^2\} \;\d{\tau}\,,
\nonumber\\
 \operatorname{Drift}(\d K) \quad&=\quad 
-\tfrac{1}{2}\min\{W^2, \alpha^2\} \;\d{\tau}\,,
\nonumber\\
 \d K\times\d H \quad&=\quad 
0\,,
\nonumber\\
 (\d H)^2 \quad&=\quad 
2 \; \d{\tau}\,,
\nonumber\\
 \operatorname{Drift}(\d H) \quad&=\quad 
-\d{\tau}\,,\label{eq:sde-K-H-planar-case-reflection}
\end{align}
where \(\d\tau=4\d t/V^2\) as before. Accordingly, for fixed \(\alpha^2\), once \(W^2\geq\alpha^2\) then \(H\) and \(K\) behave as uncorrelated Brownian motions with constant negative drifts in the \(\tau\) time-scale; moreover if \(\alpha^2>1\) then 
\[
 \operatorname{Drift}\d(\log W)/\d \tau\quad=\quad\alpha^2-1\quad
\]
 is strictly positive and so \(W\to\infty\) almost surely.
In order to argue as before we must finally show almost-sure finiteness of
\[
 T_\text{coupling}\quad=\quad
\frac{1}{4}\int_0^\infty \exp\left(2H-2K\right)\d\tau\,.
\]
Now if we scale using the ratio \((U_0/V_0)^2\) at time zero then we can
deduce the following
\undex{weak convergence}convergence-in-distribution result as \(W_0=U_0/V_0^2\to\infty\):
\begin{multline}
 \left(\frac{V_0}{U_0}\right)^2 T_\text{coupling}\quad=\quad
\frac{1}{4}\int_0^\infty \exp\left(2H-2K-(2H_0-2K_0)\right)\d\tau
\\
\quad\to\quad
\frac{1}{4}\int_0^\infty \exp\left(2\sqrt{2+\alpha^2} \, \widetilde{B}-(2-\alpha^2)\tau\right)\d\tau
\label{eq:exponential-functional}
\end{multline}
where \(\widetilde{B}\) is a standard real Brownian motion begun at \(0\). This integral is finite when \(\alpha^2<2\), so finally we deduce that coupling occurs 
in finite time for this simple mixture of reflection and synchronous coupling if we choose
\(1<\alpha^2<2\).

However we can now say much more, since the stochastic integral
in \eqref{eq:exponential-functional} is one of the celebrated and much-studied
\index{exponential functional}exponential functionals of Brownian motion
\index{Yor, M.}(\citealp{Yor-1992b}, or \citealp[p.~15]{Yor-2001}). In
particular
\index{Dufresne, D.}\citet{Dufresne-1990} has shown that such a functional
\[
 \int_0^\infty\exp\left(a\widetilde{B})_s-b s\right)\d s
\]
has the distribution of \(2/(a^2\Gamma_{2b/a^2})\), where \(\Gamma_\kappa\) is
a
\index{gamma distribution|(}Gamma-distributed random variable of index \(\kappa\). In summary,
\begin{theorem}\label{thm:reflection-synchronous-dim2}
Let \(T_\text{\rm coupling}\) be the coupling time for two-dimensional
Brownian motion plus L\'evy stochastic area under a mixture of reflection
coupling and synchronous coupling. Using the notation above,
let \(\min\{1,\alpha^2/W^2\}\) be the proportion of reflection
coupling. Scale \(T_\text{\rm coupling}\) by the square of the ratio between
initial areal difference \(U_0\) and initial spatial distance \(V_0\);
if \(1<\alpha^2<2\) then as \(W_0=U_0/V_0^2\to\infty\) so a re-scaling of the
coupling time has limiting
\index{inverse gamma distribution}Inverse Gamma distribution
\[
\left(\frac{V_0}{U_0}\right)^2 T_\text{\rm coupling}\quad\to\quad
\frac{2}{2+\alpha^2}\frac{1}{\Gamma_{(2-\alpha^2)/(2(2+\alpha^2))}}\,.
\].
\end{theorem}
Thus the
\Index{tail behaviour} of \(T_\text{coupling}\) is governed by the index of the
\index{gamma distribution|)}Gamma random variable, which here cannot exceed \(\tfrac{1}{6}\) (the limiting case when \(\alpha^2\to1\)). This index is unattainable by this means since \(H-2K\) behaves like a Brownian motion when \(\alpha^2=1\), so we cannot have \(W\to\infty\). \citet[Section 3]{Kendall-2007} exhibits a similar coupling for the planar case in which there is a state-dependent switching between reflection and synchronous coupling, depending on whether the ratio \(W\) exceeds a specified threshold; it would be interesting to calculate the tail-behaviour of the inverse of the coupling time in this case.

Note that scaling by \((U_0/V_0)^2\) rather than \(W_0=U_0/V_0^2\) quantifies
something which can be observed from detailed inspection of the
\Index{stochastic differential system} \eqref{eq:sde-V2-U2}; the rate of evolution of the areal distance \(U\) is reduced if the spatial distance \(V\) is small. The requirement \(W_0\to\infty\) is present mainly to remove the effect of higher-order terms (\(\d\tau/W^2)\) in systems such as \eqref{eq:sde-K-H-realized2}, and in particular to ensure in \eqref{eq:sde-K-H-planar-case-reflection} that \(\min\{W^2, \alpha^2\}=\alpha^2\) for all time with high probability.

In fact one can do markedly better than the reflection-synchronous mixture
coupling of Theorem \ref{thm:reflection-synchronous-dim2} by replacing
\index{coupling!synchronous coupling|)}synchronous coupling by a
\index{coupling!rotation coupling}rotation coupling for which \(\sin\theta=\sqrt{2}\,\beta/W\): similar calculations then show the index of the inverse of the limiting scaled coupling time can be increased up to the limit of \(\tfrac{1}{2}\), which remarkably is the index of the inverse of the coupling time for reflection coupling of Brownian motion alone! However this limit is not attainable by these mixture couplings, as it corresponds to a limiting case of \(\beta=2\), \(\alpha^2=3\). At this choice of parameter values \(H-2K\) again behaves like a Brownian motion so we do not have \(W\to\infty\).

In higher dimensions similar calculations can be carried out, but the geometry is more complicated; in the planar case the form of \(\underline{\underline{Z}}\) is essentially constant, whereas it will evolve stochastically in higher dimensions and relate non-trivially to \(\underline{\nu}\).
This leads to correspondingly weaker results: the inverse limiting scaled coupling time can be bounded above and below using Gamma distributions of different indices.

We should not expect these mixture couplings to be
\index{coupling!maximal coupling|(}maximal, even within the class of
\index{coupling!co-adapted coupling|(}co-adapted couplings. Indeed \citet{Kendall-2007}
gives a heuristic argument to show that maximality amongst co-adapted
couplings should be expected only when one Brownian differential is a
(state-dependent) rotation or rotated
\index{coupling!reflection coupling|)}reflection of the other. The interest of these mixture couplings lies in the ease with which one may derive limiting distributions for them, hence gaining a good perspective on how rapidly one may couple the stochastic area.

\section{Conclusion}\label{sec:conclusion}
After reviewing aspects of coupling theory, we have indicated an approach to
co-adapted coupling of Brownian motion and its stochastic areas, and shown how
in the planar case one can use Dufresne's formula to derive asymptotics of
coupling time distributions for suitable mixed couplings. Aspects of these
asymptotic distributions indicate the price that is to be paid for coupling
stochastic areas as well as the Brownian motions themselves; however it is
clear that these mixed couplings should not be expected to be
\index{coupling!maximal coupling|)}maximal amongst all co-adapted
couplings. Accordingly a very interesting direction for future research is to
develop these methods to derive estimates for
\index{coupling!time|)}coupling time distributions for more efficient
couplings using state-dependent coupling strategies as exemplified
in \citet[Section 3]{Kendall-2007}. Progress in this direction would deliver
probabilistic gradient estimates in the manner of
\index{Cranston,
M.}\citeauthor{Cranston-1991} \citetext{\citeyear{Cranston-1991}, \citeyear{Cranston-1992}}
(contrast the analytic work of
\index{Nakry, D.}\index{Baudoin, F.}\index{Bonnefont, M.}\index{Chafa\"\i, D.}\citealt{BakryBaudoinBonnefontChafai-2007}).

A further challenge is to develop these techniques for higher-\break dimensional cases. Here the two-dimensional approach extends na\"{\i}vely to deliver upper and lower bounding distributions; a more satisfactory answer with tighter bounds will require careful analysis of the evolution under the coupling of the geometry as expressed by the pair \((\nu,\underline{\underline{Z}})\).

A major piece of unfinished business in this area is to determine the extent
to which these
\index{coupling!co-adapted coupling|)}co-adapted coupling results extend to
higher-order iterated
\index{path!integral}path integrals (simple
\index{Ito, K.@It\^o, K.!It\^o calculus}It\^o calculus demonstrates that it
suffices to couple L\'evy stochastic areas in order to couple all possible
non-iterated path integrals of the form \(\int B_i\d B_j\)). Some tentative
insight is offered by the r\^ole played by the
\index{Morse, H. C.!Morse--Thue binary sequence}Morse--Thue sequence for
iterated time-integrals
\index{Price, C. J.}\citep{KendallPrice-2004}. Moreover it is possible to
generalize the invariance considerations
underlying \eqref{eq:areal-difference} for the areal difference, so as to
produce similarly invariant differences of higher-order iterated path
integrals. But at present the closing question of \citet{Kendall-2007} still
remains open, whether one can co-adaptively couple Brownian motions together
with all possible iterated path and time-integrals up to a fixed order of
iteration.\index{coupling|)}\index{Brown, R.!Brownian motion|)}\index{Levy, P.@L\'evy, P.!L\'evy stochastic area|)}

\bibliographystyle{cambridgeauthordateCMG}
 \bibliography{jfck-kendall}

\end{document}